\documentclass{amsart}
\usepackage{graphicx,amscd,color,amsmath,amsfonts,amssymb,geometry}
\newtheorem{theorem}{Theorem}[section]
\newtheorem{remark}[theorem]{Remark}

\numberwithin{equation}{section}
\begin{document}
\title[On the multilinear Bohnenblust--Hille constants: complex versus real case]{On the multilinear Bohnenblust--Hille constants: complex versus real case}
\author[Campos et al.]{J.R. Campos, D. Nu\~{n}ez-Alarc\'{o}n\textsuperscript{*}, D. Pellegrino\textsuperscript{**}, J.B. Seoane-Sep\'{u}lveda\textsuperscript{***}, D. M. Serrano-Rodr\'{\i}guez\textsuperscript{*}}

\address{Departamento de Matem\'{a}tica,\newline\indent Universidade Federal da Para\'{\i}ba,\newline\indent Rio Tinto, Brazil.}
\email{jamilsonrc@gmail.com}

\address{Departamento de Matem\'{a}tica,\newline\indent Universidade Federal da Para\'{\i}ba,\newline\indent 58.051-900 - Jo\~{a}o Pessoa, Brazil.}
\email{danielnunezal@gmail.com}
\thanks{\textsuperscript{*}Supported by Capes.}
\address{Departamento de Matem\'{a}tica,\newline\indent Universidade Federal da Para\'{\i}ba,\newline\indent 58.051-900 - Jo\~{a}o Pessoa, Brazil.}
\email{dmpellegrino@gmail.com and pellegrino@pq.cnpq.br}
\thanks{\textsuperscript{**}Supported by CNPq Edital Universal and INCT-Matem\'{a}tica.}
\address{Departamento de An\'{a}lisis Matem\'{a}tico,\newline\indent Facultad de Ciencias Matem\'{a}ticas, \newline\indent Plaza de Ciencias 3, \newline\indent Universidad Complutense de Madrid,\newline\indent Madrid, 28040, Spain.}
\email{jseoane@mat.ucm.es}
\thanks{\textsuperscript{***}Supported by grant MTM2012-34341.}
\address{Departamento de Matem\'{a}tica,\newline\indent Universidade Federal da Para\'{\i}ba,\newline\indent 58.051-900 - Jo\~{a}o Pessoa, Brazil.}
\email{dmserrano0@gmail.com}


\begin{abstract}
The results of this note arise a rupture between the behavior of the real and complex best known constants for the multilinear Bohnenblust--Hille
inequality; in one side, for real scalars, we show that new upper bounds for the real Bohnenblust--Hille inequality (the best up to now) can be obtained via a somewhat ``chaotic'' combinatorial approach, while in the complex case the combinatorial approach giving the best known constants seems to be fully controlled. We believe that the understanding of this fact is a challenging problem that may shed some new light to the subject. As a byproduct of our results we present new estimates for the constants of the Bohnenblust--Hille inequality as well as new closed formulas.
\end{abstract}

\maketitle

\section{Introduction}

Let $\mathbb{K}$ be the real or complex scalar field. The multilinear Bohnenblust--Hille inequality (\cite{Bohnenblust--Hille}, see also \cite{ba} for a more recent approach) asserts that there exists a sequence of positive scalars $\left(  B_{n}\right)  _{n=1}^{\infty}$ in $[1,\infty)$ such that
\begin{equation}
\left(  \sum\limits_{i_{1},\ldots,i_{n}=1}^{N}\left\vert U(e_{i_{^{1}}}%
,\ldots,e_{i_{n}})\right\vert ^{\frac{2n}{n+1}}\right)  ^{\frac{n+1}{2n}}\leq
B_{n}\sup_{z_{1},\ldots,z_{n}\in\mathbb{D}^{N}}\left\vert U(z_{1},\ldots
,z_{n})\right\vert \label{ul}%
\end{equation}
for all $n$-linear forms $U:\mathbb{K}^{N}\times\cdots\times\mathbb{K}%
^{N}\rightarrow\mathbb{K}$ and every positive integer $N$, where $\left(
e_{i}\right)  _{i=1}^{N}$ denotes the canonical basis of $\mathbb{K}^{N}$ and
$\mathbb{D}^{N}$ represents the open unit polydisc in $\mathbb{K}^{N}$. The
exact values for the optimal constants $B_{n}$ satisfying (\ref{ul}) remains a
mystery and are being improved throughout the time. It is worth mentioning that the
Bohnenblust--Hille inequality (and the growth of its constants) have been
shown to have applications in Quantum Information Theory (see the recent work
by Montanaro, \cite{mont}).

From now on we denote the optimal constants of the Bohnenblust-Hille inequality by $K_{n}$ (for the sake of simplicity we keep the same notation
for complex and real scalars, although the values are quite likely not the same).

The first estimates (\cite{Bohnenblust--Hille, d, Ka, Q}) suggested an exponential growth and only very recently quite different results have arisen.
The ultimate information related to the search of optimal values for constants satisfying (\ref{ul}) is:

\begin{itemize}
\item For real scalars,
\[
2^{1-\frac{1}{n}}\leq K_{n}\leq C_{n},
\]
where $\left(  C_{n}\right)  _{n=1}^{\infty}$ have a subpolynomial growth and
is given by a puzzling recursive formula (see (\ref{orig})).

\item For complex scalars,
\begin{equation}
1\leq K_{n}\leq\widetilde{C}_{n}, \label{kki}%
\end{equation}
where $\left(  \widetilde{C}_{n}\right)  _{n=1}^{\infty}$ have a subpolynomial
growth and is given by a similar recursive formula (\ref{ccxx}).
\end{itemize}

Up to now, it was an open problem, for real scalars, if $K_{n}=2^{1-\frac{1}{n}}$ or
$K_{n}=C_{n}$ or whether $K_{n}$ lies strictly between these bounds. The only
known precise value appears in the case $n=2$, since $2^{1-\frac{1}{2}}%
=K_{2}.$ For the complex case, the similar question is unsolved for the
estimates (\ref{kki}).

One of the main goals of this note is to show a somewhat surprising rupture
between the complex and real constants of the Bohnenblust--Hille inequality.
We show that the constants $C_{n}$ are not the optimal ones; in fact improved constants can be obtained via a somewhat chaotic combinatorial induction. Moreover, this result seems to be in strong contrast to the case of complex scalars, in which the best known constants are obtained via a quite controlled (and by now standard) approach.

We also provide better closed formulas for the cases of real and complex
scalars. More precisely, we show that, for all $n\geq2$,
\[
K_{n}\leq\sqrt{2}\left(  n-1\right)  ^{\log_{2}\left(  \frac{e^{1-\frac
{\gamma}{2}}}{\sqrt{2}}\right)  }\leq\sqrt{2}\left(  n-1\right)  ^{0.526322}
\]
for the case of real scalars, and
\[
K_{n}\leq\frac{2}{\sqrt{\pi}}\left(  n-1\right)  ^{\log_{2}\left(  e^{\frac
{1}{2}-\frac{1}{2}\gamma}\right)  }\leq\frac{2}{\sqrt{\pi}}\left(  n-1\right)
^{0.304975},
\]
for the case of complex scalars. Above, $\gamma \approx 0.5772$ denotes the famous Euler--Mascheroni constant.

\section{Background: the best known formulas up to today}

Let
\begin{equation}
A_{p}:=\sqrt{2}\left(  \frac{\Gamma(\frac{p+1}{2})}{\sqrt{\pi}}\right)
^{1/p}, \label{kkklll}%
\end{equation}
for $p>p_{0}\approx1.847$ and
\begin{equation}
A_{p}:=2^{\frac{1}{2}-\frac{1}{p}} \label{nob}%
\end{equation}
for $p\leq p_{0}\approx1.847.$ The exact definition of $p_{0}$ is given by the
following equality: $p_{0}\in(1,2)$ is the unique real number with
\[
\Gamma\left(  \frac{p_{0}+1}{2}\right)  =\frac{\sqrt{\pi}}{2}.
\]
The constants $A_{p}$ are the best constants satisfying Khinchine's inequality (due to Haagerup, \cite{Ha}). Up to today, the best constants satisfying the multilinear Bohnenblust--Hille inequality for real scalars appeared in \cite{jmaa} and obey the following recursive formula:
\begin{equation}
C_{m}=\left\{
\begin{array}
[c]{ll}%
1 & \text{if }m=1,\\
\left(  A_{\frac{2m}{m+2}}^{m/2}\right)  ^{-1}C_{\frac{m}{2}} & \text{if
}m\text{ is even, and }\\
\left(  A_{\frac{2m-2}{m+1}}^{\frac{-1-m}{2}}C_{\frac{m-1}{2}}\right)
^{\frac{m-1}{2m}}\left(  A_{\frac{2m+2}{m+3}}^{\frac{1-m}{2}}C_{\frac{m+1}{2}%
}\right)  ^{\frac{m+1}{2m}} & \text{if }m\text{ is odd.}%
\end{array}
\right.  \label{orig}%
\end{equation}

For complex scalars the best known constants satisfying the multilinear Bohnenblust--Hille inequality appear in \cite{jfa2}, and given by the formula

\begin{equation}
\widetilde{C}_{m}=\left\{
\begin{array}
[c]{ll}%
1 & \text{ if }m=1,\\
\left(  \left(  \widetilde{A_{\frac{2m}{m+2}}}\right)  ^{m/2}\right)
^{-1}\widetilde{C}_{\frac{m}{2}} & \text{ if }m\text{ is even, and }\\
\left(  \left(  \widetilde{A_{\frac{2m-2}{m+1}}}\right)  ^{\frac{-1-m}{2}%
}\widetilde{C}_{\frac{m-1}{2}}\right)  ^{\frac{m-1}{2m}}\left(  \left(
\widetilde{A_{\frac{2m+2}{m+3}}}\right)  ^{\frac{1-m}{2}}\widetilde{C}%
_{\frac{m+1}{2}}\right)  ^{\frac{m+1}{2m}} & \text{ if }m\text{ is odd,}%
\end{array}
\right.  \label{ccxx}%
\end{equation}
where%
\[
\widetilde{A_{p}}=\left(  \Gamma\left(  \frac{p+2}{2}\right)  \right)
^{\frac{1}{p}}.
\]

\section{New upper estimates: the exhaustive combinatorial approach}

\subsection{Real case}

Let $f:[1,2)^{2}\rightarrow\mathbb{R}$ be given by
\[
f\left(  x,y\right)  =\frac{4x-2xy}{4x+4y-4xy},
\]
$r:\mathbb{N}\rightarrow\mathbb{R}$ be defined by
\[
r\left(  x\right)  =\frac{2x}{1+x}
\]
and $A:[1,2)\rightarrow\mathbb{R}$ be given by
\[
A\left(  p\right)  =\left\{
\begin{array}[c]{c}
2^{\frac{1}{2}-\frac{1}{p}}\text{; if }p\leq p_{0}\\
\sqrt{2}\left(  \frac{\Gamma\left(  \frac{p+1}{2}\right)  }{\sqrt{\pi}
}\right)  ^{\frac{1}{p}}\text{; if }p>p_{0}.
\end{array}
\right.
\]
\bigskip

From \cite[Theorem 4.1]{defant} and using the best known constants for the
Khinchine inequality from \cite{Ha} we can see that the optimal constants
$\left(  K_{m}\right)  _{m=1}^{\infty}$ satisfying the real multilinear
Bohnenblust--Hille inequality are such that
\[
K_{m}\leq J(k,m),
\]
for all $k=1,...,\frac{m}{2}$ (when $m$ is even) and $k=1,...,\frac{m-1}{2}$
(when $m$ is odd), with
\[
J(k,m):=\left(  K_{k}\times\left(  A\left(  r\left(  k\right)  \right)
\right)  ^{k-m}\right)  ^{f\left(  r\left(  k\right)  ,r\left(  m-k\right)
\right)  }\times\left(  K_{m-k}\times\left(  A\left(  r\left(  m-k\right)
\right)  \right)  ^{-k}\right)  ^{f\left(  r\left(  m-k\right)  ,r\left(
k\right)  \right)  }.
\]
So, formally, the best estimate furnished by this method is%
\begin{equation}
\left\{
\begin{array}
[c]{c}%
K_{m}\leq P_{m}:=\min\left\{  J(k,m):k=1,...,\frac{m}{2}\right\}  \text{ if
}m\text{ is even}\\
K_{m}\leq P_{m}:=\min\left\{  J(k,m):k=1,...,\frac{m-1}{2}\right\}  \text{ if
}m\text{ is odd.}%
\end{array}
\right.  \label{fff2}%
\end{equation}
A first inspection shows that the choice%
\begin{equation}
\left\{
\begin{array}
[c]{c}%
k=\frac{m}{2}\text{ for }m\text{ even,}\\
k=\frac{m-1}{2}\text{ for }m\text{ odd}%
\end{array}
\right.  \label{cho}%
\end{equation}
seems to be the best possible (i.e., the choice where the minimum of $J(k,m)$
is achieved). For this reason, in \cite{jmaa} this approach was selected and
the formula (\ref{orig}) was presented.

As we mentioned before, at a first glance (or with the aid of some numerical
tests) it seems clear that, in general, the choice (\ref{cho}) is better than
other choices for $k$; for instance, the choice $k=2$ was investigated in
\cite{Mu}. However, in some isolated cases we now identified that this choice
of $k$ given by (\ref{cho}) was not the best one. Our main goal, rather than
just a numerical approach, is to shed light to a curious rupture between the
behavior of the best known constants for the real and complex
Bohnenblust--Hille inequalities. For this reason, in this paper we look for
the sharper constants by using the whole formula (\ref{fff2}) which comes out
with the chaotic way of generating the constants for the case of real scalars.
In view of the amount of calculations involved and since a serious precision
in the decimals is crucial, this new approach was done with a computer
program. The program, which code is in the Appendix, calculates the constants
by using the formula (\ref{fff2}). The first improvement on the constants
appear for $m=26$ and since it is a recursive procedure, this improvement
generates improvements in several other values of $m$. The following table is illustrative:%

\[%
\begin{tabular}
[c]{|l|l|l|}\hline
$m$ & {\small {new constants} $P_{m}$} & $C_{m}$ {\small {(from \cite{jmaa})}%
}\\\hline
$26$ & $<5.22772$ & $>5.22825$\\\hline
$27$ & $<5.31314$ & $>5.31447$\\\hline
$28$ & $<5.39343$ & $>5.39626$\\\hline
$29$ & $<5.47164$ & $>5.47314$\\\hline
$100$ & $<10.509$ & $>10.510$\\\hline
\end{tabular}
\ \ \ \
\]

\bigskip

From $m=27$ to $500$ the only values of $m$ for which the new constants
$P_{m}$ are not strictly smaller than $C_{m}$ are $31$, $32$, $33$, $47$,
$48$, $49$, $63$, $64$, $65$, $95$, $96$, $97$, $127$, $128$, $129$, $191$,
$192$, $193$, $255$, $256$, $257$, $383$, $384$, $385.$ As $m$ goes to
infinity, it is natural that the exhaustive combinatorial approach certainly
achieves more constants. Moreover, for certain higher values of $m,$ the
difference $C_{m}-P_{m}$ can be chosen arbitrarily large, as the following
result illustrates. The proof is simple and we omit; we just mention that the
nontrivial fact that the sequence $\left(  A_{\frac{2m}{m+2}}^{-m/2}\right)
_{m=1}^{\infty}$ is increasing (see \cite[Lemma 6.1]{jfa}) is crucial for the proof:

\begin{theorem}
Given any increasing sequence of positive real numbers $\left(  L_{j}\right)
_{j=1}^{\infty}$ with $\lim_{j\rightarrow\infty}L_{j}=\infty$, there is a
strictly increasing sequence $\left(  m_{k}\right)  _{k=1}^{\infty}$ of
positive integers so that%
\[
C_{m_{j}}-P_{m_{j}}>L_{j}\text{ for all }j.
\]

\end{theorem}

The following table illustrates the previous result:
\[%
\begin{tabular}
[c]{|c|c|}\hline
$m$ & $C_{m}-P_{m}$\\\hline
$26\cdot2^{50}$ & $>3450$\\\hline
$26\cdot2^{100}$ & $>1.19\cdot10^{11}$\\\hline
$26\cdot2^{150}$ & $>4.10\cdot10^{18}$\\\hline
\end{tabular}
\ \ \ \
\]

\subsection{Complex case: open (and deep, we believe) questions}

For complex scalars \cite[Theorem 4.1]{defant} must be replaced by a better
approach using Steinhaus variables as in \cite{jfa2}. This new approach
provides the constants $\widetilde{C}_{n}$. In this setting we made a similar
exhaustive combinatorial approach up to $m=500$ but the constants obtained
were exactly the previous from \cite{jfa2}, where the choice (\ref{cho}) is
made. Thus, it seems that no improvement can be obtained following this argument
and the question on the optimality of these constants is still open. More than
just numerical observations we feel that this rupture between the real and
complex cases is a challenging problem. Is this rupture due to some fault of
the Rademacher functions (if compared with Steinhaus variables) for the
purpose of calculating the Bohnenblust--Hille constants? What is the concrete
cause of this chaotic best choice of combinations for the real case in
contrast with an apparent perfection of the complex case?\bigskip


\section{New closed formulas}

\label{sharpest}

The notation and terminology of this section are the same as those from
\cite{jfa}, where it is proved that the optimal multilinear Bohnenblust--Hille
constants $\left(  K_{n}\right)  _{n=2}^{\infty}$ satisfy
\begin{equation}
K_{n}<1.65\left(  n-1\right)  ^{0.526322}+0.13\text{ (real scalars)}
\label{rea}%
\end{equation}
and
\[
K_{n}<1.41\left(  n-1\right)  ^{0.304975}-0.04\text{ (complex scalars).}%
\]
The proof of the above estimates is achieved by following a series of
technical steps. In the case of real scalars, using some previous lemmata, it
is observed that the sequence
\[
M_{n}=\left\{
\begin{array}
[c]{ll}%
\left(  \sqrt{2}\right)  ^{n-1} & \text{ if }n=1,2\\
DM_{\frac{n}{2}} & \text{ if }n\text{ is even, and}\\
DM_{\frac{n+1}{2}} & \text{ if }n\text{ is odd}%
\end{array}
\right.
\]
satisfies the multilinear Bohnenblust--Hille inequality, where $D=\frac
{e^{1-\frac{1}{2}\gamma}}{\sqrt{2}}.$ Then, using a \textquotedblleft uniform
approximation\textquotedblright\ argument, the estimate (\ref{rea}) is
achieved. In this section we remark that this final step of the proof, i.e.,
the uniform approximation argument, can be dropped and a quite simple argument
provides even better constants. In fact, from \cite{jfa} we know that, for all
$k\geq1$ and $n\geq2$, we have
\[
M_{n}=\sqrt{2}D^{k-1}\text{ whenever }n\in B_{k}=\{2^{k-1}+1,\ldots,2^{k}\}.
\]
Thus, $k-1\leq\log_{2}\left(  n-1\right)  $ and, hence,
\[
M_{n}\leq\sqrt{2}D^{\log_{2}\left(  n-1\right)  }=\sqrt{2}\left(  n-1\right)
^{\log_{2}\left(  \frac{e^{1-\frac{1}{2}\gamma}}{\sqrt{2}}\right)  }\leq
\sqrt{2}\left(  n-1\right)  ^{0.526322}.
\]
Using a similar argument (for complex scalars) it follows that
\[
M_{n}\leq\frac{2}{\sqrt{\pi}}\left(  n-1\right)  ^{\log_{2}\left(  e^{\frac
{1}{2}-\frac{1}{2}\gamma}\right)  }\leq\frac{2}{\sqrt{\pi}}\left(  n-1\right)
^{0.304975}%
\]
for the complex scalar field. Summarizing, we have:

\begin{theorem}
The optimal constants satisfying the Bohnenblust--Hille multilinear inequality
satisfy%
\[
K_{n}\leq\sqrt{2}\left(  n-1\right)  ^{0.526322}%
\]
for $n\geq2$ and real scalars, and%
\[
K_{n}\leq\frac{2}{\sqrt{\pi}}\left(  n-1\right)  ^{0.304975}%
\]
for $n\geq2$ and complex scalars.
\end{theorem}

Of course, the other estimates of \cite{jfa} related to the above results can
be straightforwardly improved by using these new estimates. Next, we shall
cover both real and complex cases in order to improve our previous estimates
for \emph{very large} values of $n.$

\section{Closed formulas for ``large'' values of $n$}

In this section we illustrate how the recursive essence of the best known
constants of the Bohnenblust--Hille inequality affects the calculation of
closed formulas. More precisely, we show that for big values of $n$ the
previous estimates can be pushed further. As before, the notation and
terminology of this section are the same as those from \cite{jfa}.

\subsection{Real case}

If $\left(  C_{n}\right)  _{n=1}^{\infty}$ denotes the sequence in \cite[(4.3)]{jfa}, if we fix any $k_{0},$ it is obvious that

\[
J_{n}=\left\{
\begin{array}
[c]{ll}%
C_{n} & \text{ if }n\leq2^{k_{0}},\\
DJ_{\frac{n}{2}} & \text{ if }n>2^{k_{0}}\text{ is even, and }\\
D\left(  J_{\frac{n-1}{2}}\right)  ^{\frac{n-1}{2n}}\left(  J_{\frac{n+1}{2}%
}\right)  ^{\frac{n+1}{2n}} & \text{ if }n>2^{k_{0}}\text{ is odd,}%
\end{array}
\right.
\]
with $D=\frac{e^{1-\frac{1}{2}\gamma}}{\sqrt{2}}$, satisfies the multilinear
Bohnenblust--Hille inequality. For $n>2^{k_{0}}$, let $k_{1}>k_{0}$ be such
that
\[
2^{k_{1}-1}+1\leq n\leq2^{k_{1}}.
\]
Then%
\[
k_{1}-k_{0}\leq\log_{2}\left(  \frac{n-1}{2^{k_{0}-1}}\right)
\]
and, since $(J_{n})_{n=1}^{\infty}$ is increasing, the optimal constants
$K_{n}$ satisfying the multilinear Bohnenblust--Hille inequality are so that
\begin{align*}
K_{n}  &  \leq J_{2^{k_{1}}} =D^{k_{1}-k_{0}}C_{2^{k_{0}}}\\
&  \leq C_{2^{k_{0}}}D^{\log_{2}\left(  \frac{n-1}{2^{k_{0}-1}}\right)  }\\
&  =\frac{C_{2^{k_{0}}}}{D^{k_{0}-1}}\left(  n-1\right)  ^{\log_{2}D}.
\end{align*}
We thus have
\[
K_{n}\leq\frac{C_{2^{k_{0}}}}{D^{k_{0}-1}}\left(  n-1\right)  ^{\log
_{2}\left(  \frac{e^{1-\frac{1}{2}\gamma}}{\sqrt{2}}\right)  }.
\]
From \cite[Theorem 3.1]{diana} we know that
\begin{equation}
C_{2^{k_{0}}}\leq4D^{k_{0}-4} \label{gt}%
\end{equation}
whenever $k_{0}\geq4.$ Thus,
\[
K_{n}\leq\frac{4}{\left(  \frac{e^{1-\frac{1}{2}\gamma}}{\sqrt{2}}\right)
^{3}}\left(  n-1\right)  ^{\log_{2}\left(  \frac{e^{1-\frac{1}{2}\gamma}
}{\sqrt{2}}\right)  }<1.338887\left(  n-1\right)  ^{0.526322}.
\]
Summarizing, we have:

\begin{theorem}
If $n>16$, then
\[
K_{n}\leq\frac{4}{\left(  \frac{e^{1-\frac{1}{2}\gamma}}{\sqrt{2}}\right)
^{3}}\left(  n-1\right)  ^{\log_{2}\left(  \frac{e^{1-\frac{1}{2}\gamma}
}{\sqrt{2}}\right)  }.
\]
Numerically,
\begin{equation}
K_{n}<1.338887\left(  n-1\right)  ^{0.526322}. \label{tp}%
\end{equation}

\end{theorem}

If we use the exact value of $C_{2^{k_{0}}}$ instead of estimate (\ref{gt}) we
can improve (\ref{tp}) as $n$ grows. For example,
\begin{align*}
n  &  >2^{6}\Rightarrow K_{n}<1.310883\left(  n-1\right)  ^{0.526322}\\
n  &  >2^{7}\Rightarrow K_{n}<1.306156\left(  n-1\right)  ^{0.526322}\\
n  &  >2^{8}\Rightarrow K_{n}<1.303787\left(  n-1\right)  ^{0.526322}.
\end{align*}

\subsection{Complex case}

Let $\left(  C_{n}\right)  _{n=1}^{\infty}$ denote the sequence in
\cite[Theorem 2.3]{jfa2}. If we fix any $k_{0}$, and as the authors did in
\cite{jfa}, we can show that
\[
J_{n}=\left\{
\begin{array}
[c]{ll}%
C_{n} & \text{ if }n\leq2^{k_{0}},\\
DJ_{\frac{n}{2}} & \text{ if }n>2^{k_{0}}\text{ is even, and}\\
D\left(  J_{\frac{n-1}{2}}\right)  ^{\frac{n-1}{2n}}\left(  J_{\frac{n+1}{2}%
}\right)  ^{\frac{n+1}{2n}} & \text{ if }n>2^{k_{0}}\text{ is odd,}%
\end{array}
\right.
\]
with $D=e^{\frac{1}{2}-\frac{1}{2}\gamma}$, satisfies the multilinear
Bohnenblust--Hille inequality. For $n>2^{k_{0}}$, by mimicking the real case
we obtain
\[
K_{n}\leq\frac{C_{2^{k_{0}}}}{D^{k_{0}-1}}\left(  n-1\right)  ^{\log
_{2}\left(  e^{\frac{1}{2}-\frac{1}{2}\gamma}\right)  }.
\]
Thus, using the values of $C_{2^{k}}$ from \cite{jfa2} we have
\begin{align*}
n  &  >2^{3}\Rightarrow K_{n}<1.029610\left(  n-1\right)  ^{0.304975},\\
n  &  >2^{6}\Rightarrow K_{n}<0.996322\left(  n-1\right)  ^{0.304975},\,\\
n  &  >2^{15}\Rightarrow K_{n}<0.991365\left(  n-1\right)  ^{0.304975}.
\end{align*}


\section{Open Questions}

Although the results of this note are simple to an expert, we believe that some new issues are bring into light, as the following open problems illustrate:

\begin{itemize}

\item[(1.-)] Is there any explanation for the apparent ``chaos'' in the real case in contrast with the ``perfect'' behavior in the complex case?

\item[(2.-)] Is it a fault of the Rademacher system for the purposes we need? More precisely, is there any sequence of random variables (for the real case) which behaves better (with better constants) as it happens for Steinhaus variables in the case of complex scalars?

\item[(3.-)] Are the constants obtained here the optimal ones for the Bohnenblust--Hille inequality?

\end{itemize}

\section{Appendix: the codes}

In the code below, for real scalars, note that we replaced $p_{0}$ by $1.846999$. We remark that this procedure does not cause any problem (no lack of precision in the estimates). The reason is simple. In fact, since the function $r$ is increasing and
\begin{align*}
r(12)  &  =\frac{24}{13}<1.8463<p_{0}\\
r(13)  &  =\frac{26}{14}>1.857>p_{0},
\end{align*}
there is absolutely no difference in working with $1.846999$ instead of $p_{0}.$ As a matter of fact, we could have even used $1.847$ instead of
$1.846999$. The new constants, up to $500$, can be easily checked by means of, for instance, the code given below (that was made using the
\textit{Mathematica} package and provides the first 500 values of the constants).

\begin{itemize}
\item Code for the Real case:
\small{
\begin{verbatim}
M:=500;digits=100;
f[x_,y_]:=(4*x-2*x*y)/(4*x+4*y-4*x*y);
r[x_]:=(2*x)/(1 + x);
A[k_]:=
  If[k>1.846999,
  Sqrt[2]*(Gamma[(k+1)/2]/Sqrt[Pi])^(1/k),
  2^(1/2-1/k)];
CBH[1]:=1;
CBH[2]:=Sqrt[2];
For[m=3,m<=M,m++,
 CBH[m]=
  Min[
   Table[
    N[((CBH[k]*A[r[k]]^(k-m))^(f[r[k],r[m-k]]))*
    ((CBH[m-k]*A[r[m-k]]^(-k))^(f[r[m-k],r[k]])),digits],
    {k,1,If[Mod[m,2]==0,m/2,(m-1)/2]}]
   ];
 Print[{m,CBH[m]}];
]
\end{verbatim}
}

\item Code for the Complex case:

\small{
\begin{verbatim}
M:=500;digits=100;
f[x_,y_]:=(4*x-2*x*y)/(4*x+4*y-4*x*y);
r[x_]:=(2*x)/(1+x);
A[k_]:=Gamma[(k+2)/2]^(1/k);
CBH[1]:=1;
CBH[2]:=2/Sqrt[Pi];
For[m=3,m<=M,m++,
 CBH[m]=
  Min[
  Table[
    N[((CBH[k]*A[r[k]]^(k - m))^(f[r[k],r[m-k]]))*
    ((CBH[m-k]*A[r[m-k]]^(-k))^(f[r[m-k],r[k]])),digits],
    {k, 1,If[Mod[m,2]==0,m/2,(m-1)/2]}]
  ];
 Print[{m,CBH[m]}];
]
\end{verbatim}
}
\end{itemize}

\begin{remark}
(Added in October, 7, 2013) The present paper is no longer submitted to any journal, since
the main question raised by this paper (i.e., the optimality of the constants
presented here) was recently settled (in the negative) in consequence
of a new interpolative approach to the Bohnenblust--Hille inequality,
introduced in \cite{ba}. In fact, by combining ideas of \cite{ba, q1},
Fr\'{e}d\'{e}ric Bayart (in collaboration with D. Pellegrino and J.
Seoane-Sepulveda) proved that the optimal Bohnenblust--Hille constants are so
that%
\begin{equation}
K_{n}\leq Cn^{\frac{1-\gamma}{2}}\label{87}%
\end{equation}
for complex scalars and%
\begin{equation}
K_{n}\leq Cn^{\frac{2-\gamma-\ln2}{2}}\label{90}%
\end{equation}
for real scalars and these estimates are quite better than the ones of the
present note. The proof of the new estimates is done by an adequate use of the
new interpolative procedure from \cite{ba} (it will appear in a joint paper of
Fr\'{e}d\'{e}ric Bayart, D. Pellegrino and J. Seoane-Sepulveda). The argument
is the following: from the multiple Khinchin inequality one can easily prove
that the constant associated to the Bohnenblust--Hille exponent $\left(
\frac{2n-2}{n},...,\frac{2n-2}{n},2\right)  $ is dominated by  $A_{\frac{2n-2}{n}}^{-1}K_{n-1}$. From a simple variation of Proposition
3.1 of \cite{ba}, varying the position of the power $2$ in $\left(
\frac{2n-2}{n},...,\frac{2n-2}{n},2\right)  $ we still have the upper bound
$A_{\frac{2n-2}{n}}^{-1}K_{n-1}$ for the Bohnenblust--Hille exponents $\left(
\frac{2n-2}{n},...\frac{2n-2}{n},2,\frac{2n-2}{n},...\frac{2n-2}{n}\right)  $
regardless of the position of the power $2$. Interpolating the $n$
Bohnenblust--Hille exponents
\[
\left(  \frac{2n-2}{n},...\frac{2n-2}{n},\overset{\text{position }j}{2}%
,\frac{2n-2}{n},...\frac{2n-2}{n}\right)  ,
\]
with $j=1,...,n$, with $\theta_{1}=\cdots=\theta_{n}$, one obtains the
Bohnenblust--Hille exponent $\left(  \frac{2n}{n+1},...,\frac{2n}{n+1}\right)
$ with the constant $A_{\frac{2n-2}{n}}^{-1}K_{n-1}$, i.e.,
\[
K_{n}\leq A_{\frac{2n-2}{n}}^{-1}K_{n-1}.
\]
By using the optimal estimates for $A_{\frac{2n-2}{n}}$ and properties of the
gamma function one obtains (\ref{87}) and (\ref{90}). Alternatively, using the
optimal constants of the Khinchin inequality, we have the formulas:
\end{remark}%

\[
K_{m}\leq%
{\displaystyle\prod\limits_{j=2}^{m}}
\Gamma\left(  2-\frac{1}{j}\right)  ^{\frac{-j}{2j-2}}%
\]

for complex scalars and%

\begin{equation}
K_{m}\leq2^{\frac{446381}{55440}-\frac{1}{2}m}%
{\displaystyle\prod\limits_{j=14}^{m}}
\left(  \frac{\Gamma\left(  \frac{3}{2}-\frac{1}{j}\right)  }{\sqrt{\pi}%
}\right)  ^{\frac{-j}{2j-2}}\label{11}%
\end{equation}
for real scalars (and $m\geq14$). For $2\leq m\leq13$ we have%
\[
K_{m}\leq%
{\displaystyle\prod\limits_{j=2}^{m}}
2^{\frac{1}{2j-2}}.
\]

\end{document}